\documentclass[12pt]{amsart} 

\usepackage{amssymb} 
\usepackage{amsmath} 
\usepackage{fullpage} 
\usepackage{url}

\newtheorem{theo}{Theorem}

\newtheorem{cor}{Corollary}
\theoremstyle{definition}
\newtheorem{defi}{Definition}
\theoremstyle{remark}
\newtheorem{rem}{Remark}

\def\Er{{\mathbb E}}

\def\Kr{{\mathbb K}}

\def\Pr{{\mathbb P}}
\def\Qr{{\mathbb Q}}
\def\Rr{{\mathbb R}}

\def\Ac{{\mathcal{A}}}

\def\Fc{{\mathcal{F}}}

\def\Pc{{\mathcal{P}}}

\def\Sc{{\mathcal{S}}}

\def\one{{\rm \bf 1}}
\def\esssup{\operatorname{ess.\!sup}}
\def\essinf{\operatorname{ess.\!inf}}

\def\Argmin{\operatorname{Arg\,min}}
\def\({\left(}     
\def\){\right)}    
\def\[{\left[}     
\def\]{\right]}

\def\as{{\frenchspacing a.s.}~}
\begin{document}
\title{A Remark on the Structure of Expectiles} 
\author{Freddy Delbaen}
\address{ETH Z\"urich -- D-MATH \\ R\"{a}mistrasse
   101\\ 8092 Z\"{u}rich, Switzerland}
 \address{Institut f\"ur Mathematik,
 Universit\"at Z\"urich, Winterthurerstrasse 190,
 8057 Z\"urich,
 Switzerland}
\date{First version 8 July 2013, this version \today}

\begin{abstract}
 Expectiles were defined using a minimisation principle. They form a special class of coherent risk measures.  We will describe the scenario set and we will show that there is a most severe commonotonic risk measure that is smaller than the given expectile. \end{abstract}

\maketitle

\section{Introduction}

In \cite{JZ}, Johanna F. Ziegel developed a study of statistical properties of a special class of risk measures, called expectiles.  The present paper was written after the author attended a talk by P. Embrechts on expectiles and other aspects of risk management.  The talk was held during the workshop on ``Nonlinear Expectations and Knightian Uncertainty" organised by the Institute of Mathematical Sciences of the National University of Singapore.  The author is grateful to P. Embrechts and J. Ziegel for discussions on this topic.  He also thanks the NUS for the hospitality and the many discussions during the workshop.  The author also wishes to thank E. Rosazza for pointing out that some of the results were already available in \cite{Bellini}.  The presentation here is a little bit different and some of the results are more precise.

We will freely use the theory of coherence as developed by Artzner, Delbaen, Eber and Heath, \cite{ADEH1}, \cite{ADEH2}.  For the mathematics involved we will use the monographs of the author, \cite{Pisa}, \cite{FDbook}.  To remain in the same environment and to make the bridge with classical utility theory possible, we will work with coherent utility functions.  Up to a sign change they are the same as risk measures. The model we need is essentially a one period model and all random variables will be defined on a fixed atomless probability space $(\Omega, \Fc,\Pr)$.  The hypothesis that the space is atomless is not a big restriction as this simply means that on $\Omega$ we can define e.g. a random variable that is uniformly distributed on $[0,1]$.  The space $L^\infty(\Omega,\Fc,\Pr)$ is the space of all bounded real valued random variables defined on $\Omega$ and where random variables that are equal almost surely, are identified.  Most of the time we simply write $L^\infty$.  The space of all integrable random variables (again defined modulo equality \as) is denoted by $L^1$. The duality $(L^1,L^\infty)$ plays an essential role in the theory, see \cite{FDbook}.  The reader not familiar with duality theory can find more details in books on functional analysis.
\begin{defi} A mapping $u\colon L^\infty\rightarrow \Rr$ is called a Fatou coherent utility function if
\begin{enumerate}
\item $u(0)=0$ and  $\xi\ge 0$ implies $u(\xi)\ge 0$,
\item for $a\in\Rr$ we have $u(\xi+a)=u(\xi)+a$,
\item $u(\xi+\eta)\ge u(\xi)+u(\eta)$,
\item if $\lambda \ge 0$ then $u(\lambda \xi)=\lambda u(\xi)$,
\item if $\xi_n$ is a uniformly bounded sequence tending to $\xi$ in probability, then $u(\xi)\ge\limsup u(\xi_n)$,
\item if for every  uniformly bounded sequence $\xi_n$,  tending to $\xi$ in probability, we have $u(\xi)=\lim u(\xi_n)$, we say that $u$ satisfies the Lebesgue property.
\end{enumerate}
\end{defi}
The structure of these utility functions was given in \cite{ADEH1,ADEH2} and is described in the following theorem, see \cite{FDbook}:
\begin{theo} If $u$ is Fatou and coherent then there is a a convex closed set $\Sc$ of probability measures $\Qr$, absolutely continuous with respect to $\Pr$, and such that for every $\xi\in L^\infty$:
$$
u(\xi) = \inf_{\Qr\in\Sc} \Er_\Qr[\xi].
$$
Here $\Er_\Qr$ denotes the expected value taken under the measure $\Qr$. Since $\Qr\ll\Pr$, $\Er_\Qr$ can be defined on $L^\infty(\Omega,\Fc,\Pr)$. The elements $\Qr\in\Sc$ will be identified with their Radon-Nikod\'ym derivative $\frac{d\Qr}{d\Pr}$ and hence we can write $\Sc\subset L^1$.
\end{theo}
\begin{rem}  As analysed in \cite{FDbook} the set $\Sc$ is weakly compact in $L^1$ if and only if the infimum is a minimum for every $\xi\in L^\infty$. This property is also equivalent to the Lebesgue property.
\end{rem}
Expectiles were introduced in \cite{NP} and can be defined as follows.  The operator $\Argmin$ means that we take the argument where the function attains its minimum.  To make the definition easier let us define the following function defined for $l\in\Rr$.
$$
\phi_{\tau,\xi}(l)=\tau \Er\[\(\(\xi-l\)^+\)^2\] +(1-\tau)\Er\[\(\(l-\xi\)^+\)^2\],
$$
here $0 < \tau < 1$ and $\xi\in L^\infty$. The function $\phi$ is clearly strictly convex.  So there is a unique minimiser.
\begin{defi} For given $0 <\tau < 1$, the expectile $e_\tau$ is defined as
$$
e_\tau(\xi)=\Argmin \phi_{\tau,\xi}.
$$
\end{defi}
It is known, see \cite{Bellini}, \cite{JZ} and the references therein, that for $\tau \le 1/2$, $e_\tau$ defines a Fatou coherent utility function. For $\tau = 1/2$ we find the expected value.  We leave it to the reader to analyse what happens when $\tau\rightarrow 0$. The analysis will be easier once the scenario set is known.  The following theorem is proved in  \cite{NP}.  The expression can be found using first order calculus but we emphasize that it remains valid without any hypothesis on differentiability or continuity of the distribution function of $\xi$.  The definition using $\Argmin$ needs that $\xi\in L^2$ whereas the next theorem only needs $\xi\in L^1$.
\begin{theo}  For $0<\tau\le 1/2$ and $\xi\in L^\infty$, the value $e_\tau(\xi)$ is the unique solution of the equation:
$$
\tau \Er\[\(l-\xi\)^+\]-(1-\tau)\Er\[\(l-\xi\)^-\]=0.
$$
\end{theo}
The above gives a relation with examples 2,3,4, page 34,35 in \cite{FDbook}. There it is shown that for a strictly monotone concave function $\kappa:\Rr\rightarrow \Rr$, $\kappa(0)=0$, the set $\Ac=\{\xi\mid \Er[\kappa(\xi)]\ge 0\}$ defines a concave monetary utility function on $L^\infty$. In the language of statistics these monetary utility functions are elicitable (see \cite{JZ} for details on this notion). We will not use the more general concept of concave monetary utility functions. This concept was developed in \cite{FS} and for some problems the concave case can be reduced to the case of coherent functions, \cite{Pisa}, \cite{FDbook}. The set $\Ac$ is a cone if the function $\kappa$ is positively homogeneous, i.e. is of the form $\kappa(x)= - \beta x^- + \alpha x^+$, where $0<\alpha\le\beta<\infty$.  The case $0<\tau \le 1/2 \le (1-\tau)$ fits precisely in that framework. This gives an alternative way to see that $e_\tau$ is a coherent utility function. Using $\Ac=\{\xi\mid \Er[-(1-\tau)\xi^-+\tau\xi^+]\ge 0\}$ we can define
$$
e_\tau(\xi)=\sup\{x\mid \xi-x\in\Ac\}=\max\{x\mid \xi-x\in\Ac\}.
$$
\section{The scenario set}

In \cite{FDbook}, page 35, example 4, it is shown that the scenario set of $e_\tau$ is given by
$$
\Sc=\left\{\Qr\mid \text{ there is } a>0\text{ such that }a\le\frac{d\Qr}{d\Pr}\le\frac{1-\tau}{\tau}a\right\}.
$$
Because $a\le 1$, we get that for $\Qr\in\Sc$: $\frac{d\Qr}{d\Pr}\le \frac{1-\tau}{\tau}$ and therefore the set $\Sc$ is a weakly compact set in $L^1$.  We deduce that 
$$
e_\tau(\xi)=\min\left\{\Er_\Qr[\xi]\mid \Qr\in\Sc\right\}\ge \min\left\{\Er_\Qr[\xi]\mid\frac{d\Qr}{d\Pr}\le\frac{1-\tau}{\tau}\right\}.
$$
The right hand side gives the definition of the CVaR or shortfall or tail expectation with treshold $\frac{\tau}{1-\tau}\le1$. We will give a better inequality later on. The above also shows that for an element $\Qr\in\Sc$ we must have $\frac{d\Qr}{d\Pr}\ge \frac{\tau}{1-\tau}$.  This implies that elements in $\Sc$ are bounded away from zero and hence $\Sc\neq \{\Qr \mid\frac{d\Qr}{d\Pr}\le\frac{1-\tau}{\tau}\}$.

Because the scenario set is weakly compact, the utility function $e_\tau$ has the Lebesgue property.  This means that for uniformly bounded sequences of random variables, $\xi_n$,  that converge to $\xi$ in probability, we have $\lim_n e_\tau(\xi_n)=e_\tau(\xi)$. It also shows that for every $\xi\in L^\infty$ there is an element $\Qr\in\Sc$ such that $e_\tau(\xi)=\Er_\Qr[\xi]$.
\begin{rem}  If $\tau\rightarrow 0$, then $e_\tau(\xi)\rightarrow \essinf\xi$, furthermore the mapping $\tau \rightarrow e_\tau(\xi)$ is for each $\xi\in L^\infty$, continuous and nondecreasing on $[0,1/2]$. More continuity results and differentiability properties can be found in \cite{NP}.
\end{rem}
\begin{theo} The extreme points of the set $\Sc$ are precisely the elements $\Qr$ of the form
$$
\frac{d\Qr}{d\Pr}= a\one_A + \beta a\one_{A^c},
$$
where $\Pr[A].\Pr[A^c]>0$, $\beta=\frac{1-\tau}{\tau}$ and $a=\frac{1}{\beta + \Pr[A](1-\beta)}=\frac{1}{1+(\beta-1)\Pr[A^c]}$.
\end{theo}
{\bf Proof: } The proof consists of two parts.  First we suppose that $\Qr$ is of the form described in the theorem. Suppose that $h=\frac{d\Qr}{d\Pr}=\frac{f+g}{2}$ where $f,g$ are two densities in $\Sc$.  Suppose for the moment that there is $\epsilon>0$ such that $\Pr[f\le a-\epsilon]>0$. Then we have that $\essinf f\le a-\epsilon$ and since $f\in\Sc$ we must have $\esssup f\le \beta (a-\epsilon)$.  This implies that on the set $A^c$ we must have $g\ge \beta(a+\epsilon)$ hence we must have $\essinf g\ge (a+\epsilon)$ since $g\in\Sc$.  Putting together the inequalities for $g$ we get
$$
\Er[g]\ge (a+\epsilon)\Pr[A] + \beta(a+\epsilon)\Pr[A^c] \ge a\Pr[A]+\beta a\Pr[A^c] +\epsilon\Pr[A]+\epsilon\Pr[A^c]>1.
$$
This is a contradiction and hence we get that $f\ge a$ \as, and by symmetry also  $g\ge a$ \as. This shows that $f=a$ on $A$ and hence also $g=a$ on $A$. Then we must have $f,g\le \beta a$ on $A^c$ and hence we get that $f=g=\beta a$ on $A^c$.  Summarising we have proved that $f=g=h$ or in other words $h$ is extreme.  Now suppose that $h$ is in $\Sc$ and that $a\le h \le \beta a$.  Suppose that there is $\epsilon>0$ with $\Pr[a+\epsilon\le h\le \beta(a-\epsilon)]$.  We will show  that $h$ is not extreme.  Because the probability space is atomless, there are two non negligible sets $C\cap D=\emptyset$ such that $C\cup D\subset \{a+\epsilon\le h\le \beta(a-\epsilon)\}$ and $\Pr[C]=\Pr[D]>0$.  Take $\delta=\epsilon/2$ and let $\eta=\delta\one_C-\delta\one_D$.  Clearly $h+\eta$ nd $h-\eta$ are still in $\Sc$.  This shows that $h$ is not extreme.  This shows that an extreme point is of the form $h= a\one_A + \beta a\one_{A^c}$ for some $A$.\qed
\begin{rem}  Because $1\in\Sc$ is not an extreme point (except in the trivial case $\tau=1/2$) and because $1$ is in the closure of the set of extreme points, we find that the set of extreme points is not closed in the strong (= norm) $L^1$ topology.  This is not quite uncommon in spaces of dimension strictly greater than 2 and certainly in infinite dimensional spaces.
\end{rem}
Because the extreme points allow to calculate the value of a linear functional we get
$$
e_\tau(\xi)=\inf\{\Er[\xi ( a\one_A + \beta a\one_{A^c})]\mid 0<\Pr[A]<1\}.
$$
This expression can be transformed to the situation where  $\Omega=[0,1]$.  We do not give all the details as the reader can find them in \cite{JZ}. The quantile function of a random variable has the same law as the random variable and using extreme points of $\Sc$ on $[0,1]$ that are decreasing we get:
\begin{theo} Let $q_x$ be a quantile function of $\xi$ defined for  $x\in ]0,1[$.  Then $e_\tau(\xi)$  is given by the minimum of the function
$$
x\rightarrow \int_0^1 q_u \(\frac{1}{1+(\beta-1)x}\(\one_{[x,1]}+\beta \one_{[0,x]}\)\)\,du.
$$
\end{theo}
The above can be made more explicit using first order calculus.  Since it does not give a closed form solution, we do not pursue this idea and leave it as an exercise.
\begin{theo} For $\xi =\one_C$ we get $e_\tau(\one_C)=\frac{\Pr[C]}{\beta-(\beta-1)\Pr[C]}=\frac{\Pr[C]}{\beta\Pr[C^c]+\Pr[C]}$. For $\xi =-\one_C$ we get $e_\tau(-\one_C)=\frac{-\beta\Pr[C]}{\beta-(\beta-1)\Pr[C]}=\frac{-\beta\Pr[C]}{\beta\Pr[C^c]+\Pr[C]}$.

\end{theo}
{\bf Proof } This can be shown using the previous result but it can also be done explicitly using the structure of the extreme points. So we must find the minimum of the expressions
$$\begin{aligned}
&\Er\[\one_C\frac{\one_A +\beta\one_{A^c}}{\beta-(\beta-1)\Pr[A]}\]\\
&=\frac{\Pr[A\cap C]+\beta\Pr[A^c\cap C]}{\beta-(\beta-1)\Pr[A]}\\
&=\frac{\beta\Pr[C]-(\beta-1)\Pr[A\cap C]}{\beta-(\beta-1)\Pr[A]}.
\end{aligned}
$$
We first fix the magnitude $\Pr[A]$ and try to find the optimal position of $A$.  There are two cases.  Let us start with $\Pr[A]\ge\Pr[C]$. It is clear that a set $A$ containing $C$ gives the smallest outcome. So in this case we get
$$
\frac{\Pr[C]}{\beta-(\beta-1)\Pr[A]}.
$$
Among these the value $\Pr[C]=\Pr[A]$ is the smallest, i.e. $\frac{\Pr[C]}{\beta-(\beta-1)\Pr[C]}$.  The second case is when $\Pr[A]\le\Pr[C]$.  Again we get that $A\subset C$ gives the smaller value, resulting in
$$
\frac{\beta\Pr[C]-(\beta-1)\Pr[A]}{\beta-(\beta-1)\Pr[A]}.
$$
Again using first year calculus we get that the smallest value is attained for $\Pr[A]=\Pr[C]$.

The case $-\one_C$ is handled by using the expression above for $C^c$ and using $-\one_C=\one_{C^c}-1$.\qed

\section{Commonotonicity}

In this section we will exploit the relation with commonotonicity.  The relation with shortfall was already proved above but this is not the best possible result.  Let us recall some definitions.
\begin{defi} Two elements $\xi,\eta$ in $L^\infty$ are commonotonic if for an independent copy $(\xi',\eta')$
of the couple $(\xi,\eta)$ we have $(\xi-\xi')(\eta-\eta')\ge 0$ \as. The coherent utility function $u$ is commonotonic if $u(\xi+\eta)=u(\xi)+u(\eta)$ for commonotonic $\xi,\eta$.
\end{defi}
Commonotonic coherent utility functions were characterised by David Schmeidler \cite{Schm}, and are related to convex games.  Ryff's theorem \cite{Ryff},  allows to characterise the scenario sets of commonotonic law invariant coherent utility functions. For an alternative presentation and detailed analysis see \cite{FDbook}.  Let us recall some results for law invariant functions.  This is a shortcut but the reader who wants to learn more can look up the papers by Schmeidler and by the author. The basic ingredient is a convex function (distortion) $f\colon [0,1]\rightarrow [0,1]$ such that $f(0)=0,f(1)=1$.  We will only treat the case where $f$ is continuous at $1$ since this is the case that is relevant for us. For such a distortion we introduce the set
$$
\Pc=\{\Qr\mid \forall A: \Qr[A]\ge f(\Pr[A])\}.
$$
With the function $f$ we define the Choquet integral (for $\xi\ge 0$):
$$
\int_0^\infty f(\Pr[\xi > x])\,dx = \min_{\Qr\in\Pc} \Er_\Qr[\xi].
$$
The equality between the two expressions is not immediate but follows from the general theory of convex games. Because $f$ is continuous at the point $1$, the set $\Pc$ is weakly compact in $L^1$. Ryff's theorem says that the extreme points of $\Pc$ are the elements that have the same law as $f'$ (seen as a random variable on $[0,1]$). It also follows from the general theory that all commonotonic law invariant coherent utility functions with the Lebesgue property are of this form.

To see the relation with $e_\tau$ let us look at the function $f$ defined as $f(\Pr[C])=e_\tau(\one_C)$.  This function which equals $f(x) = \frac{x}{\beta-(\beta-1)x}$ is convex and satisfies the conditions to define a utility function $v$.  This utility function is commonotonic and for indicators we have $v(\one_C)=e_\tau(\one_C)$.
\begin{theo}  With the notation above we have that for all $\xi$: $v(\xi)\le e_\tau(\xi)$, hence $\Sc \subset \Pc$.
\end{theo}
{\bf Proof }  We only have to prove the inequality for elementary functions, i.e. functions that can be written as $\xi =\sum_{k=1}^{k=n} \alpha_k \one_{A_k}$, where $\alpha_k > 0$ and $A_n\subset A_{n-1}\dots \subset A_1$. Because of superadditivity and commonotonicity of $v$ we get:
$$
e_\tau(\xi) \ge \sum_k \alpha_k e_\tau(\one_{A_k})= \sum_k \alpha_k v(\one_{A_k})=v(\xi).
$$
\begin{rem} The inequality in the theorem is also present in \cite{JZ}, see Theorem 4.1 and especially section 4.3.  The approach here is different since it is based on geometric properties of  scenario sets.
\end{rem}
\begin{cor}  The commonotonic function $v$ is the greatest commonotonic utility that is smaller than $e_\tau$.
\end{cor}
The extreme points of $\Pc$ satisfy $\frac{1}{\beta}\le f'\le \beta$ and by the Krein Milman theorem we have the same inequalities for all elements $\Qr\in \Pc$. That means that $\Pc\subset \Sc'$ where $\Sc'$ is the scenario set of an expectile defined with $\sigma$ satisfying $\frac{1-\sigma}{\sigma}=\(\frac{1-\tau}{\tau}\)^2=\beta^2$ (clearly $\sigma\le 1/2$ if $\tau\le 1/2$). As a corollary we get:
$$
e_\tau(\xi)\ge v(\xi)\ge e_\sigma(\xi).
$$
\begin{rem}  We remark that the extreme points of the set $\Pc$ as well as the extreme points of the scenario set of the tailexpectation  $\left\{\Qr\mid \frac{d\Qr}{d\Pr}\le \frac{1-\tau}{\tau}\right\}$, form a closed set in the $L^1-$topology.  This closedness property remains true for each law invariant commonotonic coherent utility function with the Lebesgue property.
\end{rem}
\begin{rem}  Because the extreme points, $h$,  of $\Pc$ satisfy $h'(0)=\frac{1}{\beta}$ and $h'(1)=\beta$, we have that $\Sc\neq\Pc$, hence $e_\tau\neq v$ (as mappings on $L^\infty$).
\end{rem}
\section{Kusuoka Representation of expectiles}
In \cite{kusuoka}, Kusuoka gave a representation of law invariant coherent measures or in our termonology law determined coherent utility functions.  The simplest of these are averages of tailexpectations, so called spectral measures.  Kusuoka has shown that every law invariant coherent utility is the infimum of over a convex set of averages of tail expectations.  More precisely there is a weak$^*$ closed convex set of probability measures on $[0,1]$, $\Sc$, such that the utility function can be represented as
$$
u(\xi)=\inf\left\{\int_{[0,1]} u_\alpha(\xi) \nu(d\alpha)\mid \nu\in\Sc \right\}.
$$
Here $u_\alpha$ is the tail expectation of treshold $\alpha$ where $u_1(\xi)=\Er[\xi]$ and $u_0(\xi)=\essinf \xi$.  To find the representation set for expectiles we proceed as in \cite{FDbook}.  We do not give the details as they are straightforward calculations.  The result is the following (recall that $\beta=\frac{1-\tau}{\tau}\ge 1$).
\begin{theo}  The expectile $e_\tau$ can be represented as
$$\begin{aligned}
e_\tau(\xi)&=\inf\left\{\int_{[0,1]} u_\alpha(\xi) \nu(d\alpha)\mid \nu\in\Sc \right\}\\
\text{ where }\Sc&=\{\nu\mid \nu\text{ a probability on } (0,1] \text{ with }\int_{(0,1]}\frac{1}{u}\,\nu(du)\le \beta \nu(\{1\}).
\end{aligned}
$$
\end{theo}
Instead of using tail expectations we can also take distortions.  Indeed every utility function of the form
$$
u(\xi)= \int_{[0,1]} u_\alpha(\xi) \nu(d\alpha),
$$
is a distortion with distortion function
$$
f(y)=\int_{[1-y,1]}\nu(d\alpha)\frac{\alpha+y-1}{\alpha}.
$$
This follows from a direct application of Fubini's theorem. Indeed for a measure $\nu$ on $(0,1]$ (so no mass at $0$) we get:
$$\begin{aligned}
\int_{(0,1]}u_\alpha(\xi)\,\nu(d\alpha)&=\int_{(0,1]}\nu(d\alpha)\frac{1}{\alpha}\int_0^\alpha q_u\,du\\
&=\int_0^1 du\,q_u\int_{[u,1]}\nu(d\alpha)\,\frac{1}{\alpha}\\
&=\int_0^1 du\,q_u \,f'(1-u)\text{ where $f$ is the distortion function defined by}\\
f(y)&=\int_0^y f'(x)\,dx\\
&=\int_{1-y}^1 f'(1-x)\,dx\\
&=\int_{1-y}^1 dx \int_{[x,1]}\nu(d\alpha)\,\frac{1}{\alpha}\\
&=\int_{[1-y,1]}\nu(d\alpha)\,\frac{1}{\alpha}\int_{1-y}^\alpha dx\\
&=\int_{[1-y,1]}\nu(d\alpha)\,\frac{\alpha+y-1}{\alpha}.
\end{aligned}
$$
One can easily check that the function $f$ so defined is indeed convex, $f\colon [0,1]\rightarrow [0,1]$, $f(0)=0,f(1)=1$.  The convexity follows from $y\rightarrow f'(y)=\int_{1-y}^1\nu(d\alpha)\frac{1}{\alpha}$ is nondecreasing.
To fit in the representation of expectiles we must translate the condition on the densities and this becomes
\begin{theo}  The expectile $e_\tau$ is represented as
$$
e_\tau(\xi)=\inf\left\{\int_0^1 du\,q_u \, f'(1-u)\mid f\in\Fc _\beta\right\},
$$
where $\Fc_\beta$ is the (convex) set of convex functions $f\colon[0,1]\rightarrow [0,1]$, $f(0)=0;f(1)=1$ and $f'(1)\le \beta f'(0)$.
\end{theo}
\begin{rem} The above shows that there is no smallest utility of spectral type that is bigger than the expectile.  Indeed it would be realised by a distortion function $\phi$, that necessarily must be in $\Fc_\beta$.  Let us denote $v_\phi$ the coherent utility defined with the distortion $\phi$. Since $e_\tau$ can be calculated as an infimum over all such distortion functions, we get that either $e_\tau(\xi)<v_\phi(\xi$ for some $\xi$ or $e_\tau=v_\phi$.  In the former case there must be another distortion, $\psi$ that satisfies $v_\psi(\xi)<v_\phi(\xi)$, showing that $v_\phi$ was not 
the smallest element.  In the latter case $e_\tau$ would be commonotone and this contradicts the remark 6 made after theorem 6.
\end{rem}
\begin{rem}  We remark that the pointwise convergence on $\Fc$ is equivalent to the uniform convergence and that $\Fc$ is compact for this topology.
\end{rem}
\section{Continuity of $e_\tau$}
 All the above inequalities were shown for elements of $L^\infty$.  Because all the scenario sets were subsets of $L^\infty$ we  easily see that the inequalities remain true for elements in $L^1$.  This allows to recover the continuity theorem with respect to the Wasserstein metric as e.g. proved in \cite{Bellini}.  Let us recall that the Wasserstein metric is defined for probability laws on metric spaces but that for our purpose we restrict to the case $\Rr$. The definition is as follows.  We take two laws $\mu,\nu$ on $\Rr$ for which there exist absolute moments of order $p$.  Then for $0<p$ we define
 $$
 d_p(\nu,\mu)=\inf\left\{\Vert \xi-\eta\Vert_p^p=\Er[|\xi-\eta |^p]\mid \xi, \eta \text{ have law }\mu,\nu \right\}.
 $$
 The random variables $\xi,\eta$ can be defined on any atomless probability space.  The case $p=1$ is of special interest. Since the expectiles only depend on the law of the random variables, we can write $e_\tau(\mu)$ when the law of $\xi$ is $\mu$.  The first Lipschitz property in the following theorem, was proved in \cite{Bellini}.
 \begin{theo}  The monetary utility function $e_\tau$ is Lipschitz continuous for the Wasserstein metric and the Lipschitz constant is $\beta$.  On a set where all random variables have the same mean, the Lipschitz constant can be reduced to $(\beta-1)/2$.
 \end{theo}
{\bf Proof }  The proof of this  is straightforward using the scenario set. Indeed as easily seen for random variables in $L^1$ we have $|e_\tau(\xi)-e_\tau(\eta)| \le \beta \Vert\xi-\eta\Vert_1$
$$
|e_\tau(\xi)-e_\tau(\eta)| \le \sup_{\Qr\in\Sc}\Er_\Qr[ |\xi-\eta|]\le\beta \Vert\xi-\eta\Vert_1.
$$
To get the better constant we proceed as follows. Let us suppose that $\Er[\xi]=\Er[\eta]$. Then we write
$$
\begin{aligned}
e_\tau(\xi)-e_\tau(\eta)&=\min_{\Qr\in\Sc}\Er_\Qr[\xi]-\min_{\Qr\in\Sc}\Er_\Qr[\eta]\\
&=\min_{\Qr\in\Sc}\Er_\Qr[\xi] - \Er_{\Kr_0}[\eta]\text{ for a well chosen } \Kr_0\in\Sc\\
&\le \Er_{\Kr_0}[\xi]-\Er_{\Kr_0}[\eta]\\
&\le \max_{\Qr\in\Sc, \text{extreme}}\Er_{\Qr}[\xi-\eta].
\end{aligned}
$$
We will now find a bound for the last expression.  Recall that an extreme point has a density equal to
$$
\frac{\one_A+\beta\one_{A^c}}{\beta\Pr[A^c]+\Pr[A]}.
$$
Let $k=\frac{\beta+1}{2(\beta\,\Pr[A^c]+\Pr[A])} $.  Since $\Er[\xi-\eta]=0$ we also get
$$\begin{aligned}
\Er_{\Qr}[\xi-\eta]&=\Er\[\(\frac{d\Qr}{d\Pr}-k\)(\xi-\eta)\]\\
&\le \left\Vert \frac{d\Qr}{d\Pr}-k \right\Vert_\infty \,\Vert \xi-\eta\Vert_1\\
&\le \frac{\beta-1}{2(\beta\,\Pr[A^c]+\Pr[A])}\,\Vert \xi-\eta\Vert_1\\
&\le \(\frac{\beta-1}{2}\)\Vert \xi-\eta\Vert_1.
\end{aligned}
$$
The inequality in the other direction is proved by symmetry.\qed
\begin{rem}  By taking $\eta=0$, i.e. $\nu=\delta_0$ and $\xi$ of the form $-\one_C$, we can see that the constant $\beta$ is optimal.   Since the second inequality with the constant $\frac{\beta-1}{2}$ implies (by subtracting and adding $\Er[\xi-\eta]$),  the first inequality with constant $\beta$ and since the latter constant is best possible, we get that also $\frac{\beta-1}{2}$ must be best possible.
\end{rem}
\begin{rem} We can also use the duality $(L^p,L^q)$ and get an inequality with respect to the Wasserstein metric $d_p$. The only thing to be done is to calculate the upper bound for $\Vert \frac{\one_A+\beta\one_{A^c}}{\beta -(\beta-1)\Pr[A]}\Vert_q$, a {\it straightforward } exercise.  And for the amateurs:  one can also use a duality between Orlicz spaces $(L^\Phi,L^\Psi)$ to get similar inequalities.
\end{rem}
\section{The Concave Case}
As already mentioned above we can start with an increasing concave utility functon $\phi:\Rr\rightarrow\Rr$, $\phi(0)=0$.  Then we define the acceptance set $\Ac=\{\xi\mid \Er[\phi(\xi)]\ge 0\}$.  The penalty function related to $\Ac$ can then be calculated using a variational argument.  Unfortunately there are no closed form solutions except in easy cases.  For instance we get for $\phi(x)=1-\exp(-x)$ that the penalty function is $c(\Qr)=\Er\[\frac{d\Qr}{d\Pr}\log\(\frac{d\Qr}{d\Pr}\)\]$.  Important is that these kind of acceptance sets also leads to elicitable utility functions.  We do not proceed this discussion.

\end{document}